\documentclass[preprint,12pt]{elsarticle}
\usepackage{amssymb,amsthm,amsmath,verbatim,enumerate,ifthen}
\usepackage[mathscr]{eucal}
\usepackage[utf8]{inputenc}
\usepackage[T1]{fontenc}

\def\N{\mathbb{N}}
\def\R{\mathbb{R}}

\def\C{\mathscr{C}}
\def\E{\mathscr{E}}
\def\G{\mathscr{G}}
\def\K{\mathscr{K}}
\def\L{\mathscr{L}}

\def\conv{\mathop{\mbox{\rm conv}}\nolimits}
\def\ext{\mathop{\mbox{\rm ext}}\nolimits}

\def\intr{\mathop{\mbox{\rm int}}\nolimits}

\def\id{\mathop{\mbox{\rm id}}\nolimits}
\def\ri{\mathop{\mbox{\rm ri}}\nolimits}

\newtheorem{theorem}{Theorem}[section]
\long\def\Thm#1#2{\ifthenelse{\equal{#1}{*}}{\begin{theorem*}#2\end{theorem*}}
             {\begin{theorem}\label{T#1}#2\end{theorem}}}
\def\thm#1{Theorem~\ref{T#1}}

\newtheorem{proposition}[theorem]{Proposition}
\long\def\Prp#1#2{\ifthenelse{\equal{#1}{*}}{\begin{proposition*}#2\end{proposition*}}
             {\begin{proposition}\label{P#1}#2\end{proposition}}}

\newtheorem{corollary}[theorem]{Corollary}
\long\def\Cor#1#2{\ifthenelse{\equal{#1}{*}}{\begin{corollary*}#2\end{corollary*}}
             {\begin{corollary}\label{C#1}#2\end{corollary}}}
\def\cor#1{Corollary~\ref{C#1}}

\newtheorem{lemma}[theorem]{Lemma}
\long\def\Lem#1#2{\ifthenelse{\equal{#1}{*}}{\begin{lemma*}#2\end{lemma*}}
             {\begin{lemma}\label{L#1}#2\end{lemma}}}
\def\lem#1{Lemma~\ref{L#1}}

\theoremstyle{definition}
\newtheorem{definition}[theorem]{Definition}
\long\def\Defn#1#2{\ifthenelse{\equal{#1}{*}}{\begin{definition*}\rm #2\end{definition*}}
             {\begin{definition}\label{D#1}\rm #2\end{definition}}}

\newtheorem{remark}[theorem]{Remark}
\long\def\Rem#1#2{\ifthenelse{\equal{#1}{*}}{\begin{remark*}\rm #2\end{remark*}}
             {\begin{remark}\label{R#1}\rm #2\end{remark}}}

\newtheorem{example}{Example}
\long\def\Exa#1#2{\ifthenelse{\equal{#1}{*}}{\begin{example*}\rm #2\end{example*}}
             {\begin{example}\label{Ex#1}\rm #2\end{example}}}

\def\eq#1{{\rm(\ref{E#1})}}
\def\Eq#1#2{\ifthenelse{\equal{#1}{*}}
  {\begin{equation*}\begin{aligned}#2\end{aligned}\end{equation*}}
  {\begin{equation}\begin{aligned}\label{E#1}#2\end{aligned}\end{equation}}}

\def\comment#1{}

\begin{document}
\begin{frontmatter}

\title{Support theorems in abstract settings}
\author{Andrzej Olbry\'s}
\address{Institute of Mathematics, University of Silesia,
40-007 Katowice, ul. Bankowa 14, Poland}
\ead{andrzej.olbrys@math.us.edu.pl}

\author{Zsolt P\'ales\footnote[1]{This research has been supported by the Hungarian Scientific Research Fund (OTKA) 
Grant K111651.}}
\address{Institute of Mathematics, University of Debrecen,
H-4032 Debrecen, Egyetem t\'er 1, Hungary}
\ead{pales@science.unideb.hu}


\begin{abstract}
In this paper we establish a general framework in which the verification of support theorems for generalized convex 
functions acting between an algebraic structure and an ordered algebraic structure is still possible. As for the 
domain space, we allow algebraic structures equipped with families of algebraic operations whose 
operations are mutually distributive with respect to each other. We introduce several new concepts in such 
algebraic structures, the notions of convex set, extreme set, and interior point with respect to a given family of 
operations, furthermore, we describe their most basic and required properties. In the context of the range space, we 
introduce the notion of completeness of a partially ordered set with respect to the existence of the infimum of lower 
bounded chains, we also offer several sufficient condition which imply this property. For instance, the order generated 
by a sharp cone in a vector space turns out to possess this completeness property. By taking several particular cases, 
we deduce support and extension theorems in various classical and important settings.
\end{abstract}
\end{frontmatter}

\section{Introduction}

Support theorems play crucial roles in many branches of analysis, algebra and geometry. Roughly speaking, such theorems 
lead to the representation of convex functions as the pointwise maximum of affine functions, subadditive functions as 
the pointwise maximum of additive functions, convex sets as the intersection of half spaces. The nonemptyness of the 
subgradient of a convex function at a given point (in the sense of convex analysis) can also be obtained by using a 
certain support theorem. A typical method to prove support theorems is to use the Hahn--Banach extension theorem or 
sandwich theorem or one of their generalizations to the setting of groups or semigroups (see \cite{Bad06b}, 
\cite{Bad06a}, \cite{Bal01}, \cite{Fuc74}, \cite{FucLus81}, \cite{GajSmaSma92}, \cite{JarPal15}, 
\cite{Kau66}, \cite{Kuc85}, \cite{Olb11b}, \cite{Olb17a}, \cite{Olb17b}, \cite{Pal89a}, \cite{Pal89b}, \cite{Pal89d}, 
\cite{Pal01b}, \cite{RodBou74}). A survey on these developments was given by Buskes \cite{Bus93}. The celebrated 
sandwich theorem of Rod\'e \cite{Rod78}, the abstract extension of the Hahn--Banach theorem to setting of convexity 
defined in terms of families of commuting operations, is still one of the most powerful tools. There have been many 
attempts to simplify its proof, to generalize its content and to find valuable applications (see \cite{ChaVol91}, 
\cite{FucKon80}, \cite{Kon80b}, \cite{Kon87b}, \cite{Kuh84}, \cite{Kuh87}, \cite{Pal98d}, \cite{VolWei81}). 

In the extensions and generalizations of the classical Hahn--Banach theorems, the algebraic structure of the domain 
basically did not cause any problem, sandwich theorems for extended real-valued functions over algebraic structures 
with many operations have been established. In the case of functions with values in ordered vector spaces, 
Rodrigues--Salinas and Bou \cite{RodBou74} showed that sandwich type results can only be expected for ordered vector 
spaces where the intervals have the so-called binary intersection property. Generalizations of the Hahn--Banach 
extension theorem in many settings can be deduced from sandwich theorems, however, they can be extended to operators 
with values in vector spaces with the least upper bound property, one of such an extensions is known as the 
Hahn--Banach--Kantorovi\'c theorem (see \cite{Jam70}, \cite{LusTho83}). As it was proved by Silverman and Yen 
\cite{SilYen59} (see also \cite{AngLem74a}, \cite{BonSil66}, \cite{BonSil67}, \cite{Iof81}, \cite{Nik91b}, \cite{To70}, 
\cite{To71}) the least upper bound property of the range space is indispensable, more precisely, an ordered vector space 
has the Hahn--Banach extension property if and only if it possesses the least upper bound property.

As support theorems until now have been deduced from sandwich type theorems or from Hahn--Banach type extension 
theorems, they did exist only for extended real-valued functions or vector-valued functions mapping into a space with 
the least upper bound property. In a recent paper of the first author \cite{Olb15d}, a support theorem was 
found for the vector-valued setting, namely for delta $(s,t)$-convex mappings. It turns out that delta 
$(s,t)$-convexity can be reformulated as a convexity property with respect to the Lorenz cone. However, the order 
induced by the Lorenz cone typically does not fulfills the  least upper bound property. Therefore, it turned out that 
support theorems may be obtained under much weaker conditions concerning the range space.

The main goal of this paper is to establish a general framework in which the verification of support theorems is still 
possible. As for the domain space, we allow algebraic structures equipped with families of algebraic operations whose 
operations are mutually distributive with respect to each other. (This property is much more general than the
pairwise commutativity which was needed for the setting of the Rod\'e Theorem.) We introduce several new concepts in 
such algebraic structures, the notions of convex set, extreme set, and interior point with respect to a given family of 
operations, furthermore, we describe their most basic and required properties. We  mention that no topological 
assumptions are needed, the usual conditions related the topological interior of the domain are replaced by a new 
intrinsic notion which is purely derived from the given algebraic operations. In the context of the range space, we 
introduce the notion of completeness of partially ordered set with respect to the existence of the infimum of lower 
bounded chains (which is much weaker than the existence of the infimum of lower bounded sets), we also offer several 
sufficient condition which imply this property. For instance, the order generated by a sharp cone in a vector space turns 
out to possess this completeness property.

\section{Convexity and Extremality with Respect to Families of Algebraic Operations}

The notions that we introduce below are intuitively motivated by the standard concepts that are widely used 
and applied in the theory of convex sets. This will be made transparent when we consider various particular 
cases of our definitions in the sequel.

In order to introduce the general definition of convex and extreme sets, let $\Gamma$ denote a nonempty set 
and let $n:\Gamma\to\N$ be a (so-called \emph{arity}) function throughout the rest of this paper. 

For a nonempty set $X$ and for a given family of operations on $X$
\Eq{om}{
 \omega=\big\{\omega_\gamma:X^{n(\gamma)}\to X\mid\gamma\in\Gamma\big\},
}
we say that $E\subseteq X$ is \emph{$\omega$-convex} if  
\Eq{oc}{
   \omega_\gamma(E^{n(\gamma)})\subseteq E\qquad (\gamma\in\Gamma).
}
Another notion that will play a key role in our investigations is the concept of an extreme set. We say that a 
subset $E\subseteq X$ is \emph{$\omega$-extreme} if 
\Eq{oe}{
   \omega^{-1}_\gamma(E)\subseteq E^{n(\gamma)}\qquad (\gamma\in\Gamma).
}
A point $p\in X$ is said to be \emph{$\omega$-extreme} if the singleton $\{p\}$ is an $\omega$-extreme set.
Trivially, the entire set $X$ and the empty set are $\omega$-convex and $\omega$-extreme sets. The 
collection of all $\omega$-convex subsets and $\omega$-extreme subsets of $X$ will be denoted by 
$\C_\omega(X)$ and $\E_\omega(X)$, respectively. 

We have the following easy-to-prove result.

\Prp{ce}{Let $\omega$ be a family of operations given by \eq{om}. Then $\C_\omega(X)$ is closed under 
the intersection (resp.\ under the union) of arbitrary collections (resp.\ chains) of subsets of $X$ and 
$\E_\omega(X)$ is closed under the intersection and under the union of arbitrary collections of subsets of 
$X$.}

This proposition allows us to set the following definition: The \emph{$\omega$-convex hull} $\conv_\omega(H)$ 
of a set $H\subseteq X$ is the intersection of all $\omega$-convex sets containing $H$, 
that is, $\conv_\omega(H)$ is the smallest $\omega$-convex set including the set $H$:
\Eq{*}{
\conv_{\omega}(H):=\bigcap\{E\in\C_\omega(X) \mid H\subseteq E\}.
}
Analogously, for a given set $H\subseteq X$, we may define the set $\ext_{\omega}(H)$, the 
\emph{$\omega$-extreme hull} of $H$, as the smallest (with respect to inclusion) $\omega$-extreme set 
containing $H$. In other words, 
\Eq{*}{
\ext_{\omega}(H):=\bigcap\{E\in\E_\omega(X) \mid H\subseteq E\}.
}

The following assertion easily follows from the definitions of $\omega$-convexity and $\omega$-extremality.

\Prp{05}{For arbitrary $H, H_{1}, H_{2}\subseteq X$ and sets of operations $\omega$ 
the following properties are satisfied: 
\begin{enumerate}[(1)]
 \item $\conv_{\omega}(H)=\conv_{\omega}(\conv_{\omega}(H))$ \ and \ 
 $\ext_{\omega}(H)=\ext_{\omega}(\ext_{\omega}(H))$,
 \item If $H_{1}\subseteq H_{2}$ \ then \ $\conv_{\omega}(H_{1})\subseteq \conv_{\omega}(H_{2})$ \ and \ 
 $\ext_{\omega}(H_{1})\subseteq \ext_{\omega}(H_{2})$
 \item $\conv_{\omega}(H_{1})\cup \conv_{\omega}(H_{2})\subseteq \conv_{\omega}(H_{1}\cup H_{2})$ \ and \\ 
 $\ext_{\omega}(H_{1})\cup \ext_{\omega}(H_{2})\subseteq \ext_{\omega}(H_{1}\cup H_{2})$,
 \item $\conv_{\omega}(H_{1}\cap H_{2})\subseteq \conv_{\omega}(H_{1})\cap \conv_{\omega}(H_{2})$ \ and \\
 $\ext_{\omega}(H_{1}\cap H_{2})\subseteq \ext_{\omega}(H_{1})\cap \ext_{\omega}(H_{2})$.
\end{enumerate}}

For computing the $\omega$-convex and $\omega$-extreme hulls of a set, the following result can be useful.

\Thm{hf}{Let $\omega$ be a family of operations given by \eq{om}. Then, for any subset $H\subseteq X$, we have
that
\Eq{hf1}{
  \conv_{\omega}(H)=\bigcup_{k=0}^\infty C_k \qquad\text{and}\qquad \ext_{\omega}(H)=\bigcup_{k=0}^\infty D_k,
}
where the sequences $(C_k)$ and $(D_k)$ are defined by the following recursions:
\Eq{*}{
   C_0:=H,\quad C_{k+1}&:=C_k\cup\bigg(\bigcup_{\gamma\in\Gamma} 
                  \omega_\gamma\big(C_k^{n(\gamma)}\big)\bigg),\\
   D_0:=H,\quad D_{k+1}&:=D_k\cup\bigg(\bigcup_{\gamma\in\Gamma} 
      \big\{\{x_1,\dots,x_{n(\gamma)}\}\mid \omega_\gamma(x_1,\dots,x_{n(\gamma)})\in D_k\big\}\bigg).
}}

\begin{proof} First, we prove by induction on $k$, that 
\Eq{cd}{
   C_k\subseteq \conv_{\omega}(H) \qquad\text{and}\qquad D_k\subseteq \ext_{\omega}(H),
}
which will show that both relations in \eq{hf1} hold with the inclusion ``$\supseteq$''. These statements are 
obvious for $k=0$. Assume that \eq{cd} is valid for some $k$. 

If $x\in C_{k+1}\setminus C_k$ then there exists $\gamma\in\Gamma$ and $x_1,\dots,x_{\gamma(n)}\in 
C_k\subseteq\conv_{\omega}(H)$ such that $x=\omega_\gamma(x_1,\dots,x_{n(\gamma)})$. The set 
$\conv_{\omega}(H)$ being $\omega$-convex, we have that $\omega_\gamma(x_1,\dots,x_{n(\gamma)})\in 
\conv_{\omega}(H)$, showing that $x\in \conv_{\omega}(H)$. Thus, we have obtained that
$C_{k+1}\subseteq\conv_{\omega}(H)$.
 
Now let $x\in D_{k+1}\setminus D_k$. Then there exists $\gamma\in\Gamma$ and $x_1,\dots,x_{\gamma(n)}\in X$ 
such that $x\in\{x_1,\dots,x_{n(\gamma)}\}$ and $\omega_\gamma(x_1,\dots,x_{n(\gamma)})\in 
D_k\subseteq\ext_{\omega}(H)$. The set $\ext_{\omega}(H)$ is $\omega$-extreme, hence 
$\{x_1,\dots,x_{n(\gamma)}\}\subseteq\ext_{\omega}(H)$, which implies that $x\in \ext_{\omega}(H)$. Thus, we 
have verified that $D_{k+1}\subseteq\ext_{\omega}(H)$.

For the proof of the reversed inclusions in \eq{hf1}, it suffices to show that the right hand sides of these 
relations, denoted by $C$ and $D$, are $\omega$-convex and $\omega$-extreme sets that contain $H$, 
respectively. The property that these sets contain $H$ is trivial since $H=C_0=D_0$. 

Let $\gamma\in\Gamma$ and let $x_1,\dots,x_{\gamma(n)}\in C=\bigcup_{k=0}^\infty C_k$. Then, there exists 
$k_0$ such that $x_1,\dots,x_{\gamma(n)}\in C_{k_0}$. Therefore,
\Eq{*}{
  \omega_\gamma(x_1,\dots,x_{n(\gamma)})
    \in\omega_\gamma\big(C_{k_0}^{n(\gamma)}\big)\subseteq C_{k_0+1}\subseteq C.
}
This completes the proof of the $\omega$-convexity of $C$.

To show the $\omega$-extremality of the set $D$, let $\gamma\in\Gamma$ and let $x_1,\dots,x_{\gamma(n)}\in X$ 
such that $\omega_\gamma(x_1,\dots,x_{n(\gamma)})$ belongs to $D$. Then, there exists $k_0$ such that 
$\omega_\gamma(x_1,\dots,x_{n(\gamma)})\in D_{k_0}$. By the the construction of the sequence $(D_k)$, this 
yields that $\{x_1,\dots,x_{n(\gamma)}\}\subseteq D_{k_0+1}\subseteq D$, whence the $\omega$-extremality of 
$D$ follows.
\end{proof}

The next proposition shows that the complements of $\omega$-extreme sets behave like ideals with respect to 
the operations of the family $\omega$.

\Thm{03}{Let $\omega$ be a family operations given by \eq{om}. If $E \subseteq X$ is an $\omega$-extreme set 
then, for all $\gamma\in\Gamma$ and for all $i\in\{1,\dots,n(\gamma)\}$, 
\Eq{oi}{
\omega_\gamma\big((X\setminus E)^{n(\gamma)}\big)
\subseteq
\omega_\gamma\big(\{(x_1,\dots,x_{n(\gamma)})\in X^{n(\gamma)}\mid x_i\in X\setminus E\}\big)
\subseteq X\setminus E.
}
As a consequence, $X\setminus E$ is $\omega$-convex.}

\begin{proof}
Let $\gamma\in\Gamma$ and $i\in\{1,\dots,n(\gamma)\}$. The left hand side of inclusion in \eq{oi} is trivial. If 
the right hand side inclusion in \eq{oi} were not valid, then, for some elements $x_1,\dots,x_{n(\gamma)}\in 
X$ with $x_i\in X\setminus E$, we have that $\omega_\gamma(x_1,\dots,x_{n(\gamma)})\not\in X\setminus E$, 
i.e., $\omega_\gamma(x_1,\dots,x_{n(\gamma)})\in E$. In view of the extremality of $E$, this implies that 
$(x_1,\dots,x_{n(\gamma)})\in E^{n(\gamma)}$, which contradicts $x_i\not\in E$.
\end{proof}

Now, we define a counterpart of the notion of the relative interior in terms of $\omega$-extreme points. A 
point $p\in X$ is said to be \emph{$\omega$-internal} if $\ext_{\omega}(\{p\})= X.$ The set of 
\emph{$\omega$-internal} points of $X$ is called the \emph{$\omega$-interior of $X$} and is denoted by 
$\intr_{\omega}(X),$ that is,
\Eq{*}{
\intr_{\omega}(X):=\{p\in X\mid\ext_{\omega}(p)= X\}.
}
The complement of the $\omega$-interior of $X$ is termed the \emph{$\omega$-boundary of $X$} and is 
denoted by $\partial_{\omega}(X)$.

\Prp{06}{Let $\omega$ be a family of operations given by \eq{om} and assume that 
$\intr_\omega(X)\not=\emptyset$. Then the set $\partial_{\omega}(X)$ is the largest, proper $\omega$-extreme 
subset of $X.$}

\begin{proof} By the assumption $\intr_\omega(X)\not=\emptyset$, we have that $\partial_{\omega}(X)$ is 
a proper subset of $X$. First we prove that $\partial_{\omega}(X)$ is an $\omega$-extreme subset of $X.$
Let $\gamma\in\Gamma$, let $x_{1},\dots,x_{n(\gamma)} \in X$ and assume that 
$\omega_\gamma(x_{1},\dots,x_{n(\gamma)})\in \partial_{\omega}(X)$, that is,  
$\omega_\gamma(x_{1},\dots,x_{n(\gamma)})$ is not an $\omega$-internal point of $X.$ Then the set 
$E:=\ext_{\omega}(\{\omega_\gamma(x_{1},\dots,x_{n(\gamma)})\})$ is a proper subset of $X$. 
By its $\omega$-extremality, $E$ must contain $x_i$ for all $i\in\{1,\dots,n(\gamma)\}$.
Therefore,
\Eq{*}{
\ext_{\omega}(\{x_{i}\})\subseteq \ext_{\omega}(\{\omega_\gamma(x_{1},\dots,x_{n(\gamma)})\})=E\subsetneq X,
}
which shows that $\ext_{\omega}(\{x_{i}\})$ is also a proper subset of $X$. This completes the proof of the 
inclusions $x_{i}\in\partial_{\omega}(X)$ for all $i\in\{1,\dots,n(\gamma)\}$, whence the $\omega$-extremality
of $\partial_{\omega}(X)$ follows. 

Let $F$ be a proper $\omega$-extreme subset of $X$. If $x\in F$ then $\ext_{\omega}(\{x\})\subseteq 
F\subsetneq X.$ This yields that $x\notin \intr_{\omega}(X)$, that is, $x\in \partial_{\omega}(X).$ Hence 
$F\subseteq \partial_{\omega}(X).$
\end{proof}

\comment{\Prp{01}{For an operation $\omega:X^n\to X$ and any $x_{1},\dots,x_{n}\in X$ we have 
\Eq{*}{
\ext_{\omega}(\{x_{1},\dots,x_{n}\})\subseteq \ext_{\omega}(\{\omega(x_{1},\dots,x_{n})\}).
}}
\begin{proof} Since $\omega(x_{1},\dots,x_{n})\in \ext_{\omega}(\{\omega(x_{1},\dots,x_{n})\})$ then 
$x_{1},\dots,x_{n}\in 
\ext_{\omega}(\{\omega(x_{1},\dots,x_{n})\})$, therefore, 
\Eq{*}{
\ext_{\omega}(\{x_{j}\})\subseteq  \ext_{\omega}(\{\omega(x_{1},\dots,x_{n})\}),\quad j=1,\dots,n,
}
hence
\Eq{*}{
\bigcup_{j=1}^{n}\ext_{\omega}(\{x_{j}\})\subseteq \ext_{\omega}(\{\omega(x_{1},\dots,x_{n})\}).
}
It follows from the previous proposition that $\bigcup_{j=1}^{n}\ext_{\omega}(\{x_{j}\})$ is an extreme set, 
therefore 
\Eq{*}{
\ext_{\omega}(\{x_{1},\dots,x_{n}\}) \subseteq \bigcup_{j=1}^{n}\ext_{\omega}(\{x_{j}\}).
}
\end{proof}}

\Exa{01}{
Let $X=[0,1]$ and let $\omega : [0,1]^2 \to [0,1]$ be given by formula $\omega(x,y):=\frac{x+y}{2}$. Then 
$\{0\}$, $\{1\}$, $\{0,1\}$, and $[0,1]$ are the only $\omega$-extreme sets.
}

Indeed, if $p\not\in\{0,1\}$ and $E$ is an $\omega$-extreme set containing $p$, then, for arbitrary $x\in [0,1]$, 
say $p<x$, we can choose a natural number $n$ such that $\frac{1}{n}<\min\{p,1-p\}$. Then consider the 
following sequence of points 
\Eq{*}{
x_{j}:=p+j \frac{x-p}{n}, \quad j=-1,0,1,\dots,n.
}
Since $p=\frac{x_{-1}+x_{1}}{2}$, therefore, by the $\omega$-extremality of $E$, we have that $x_{-1}, x_{1} \in E$. 
Analogously, because $x_{1}=\frac{x_{0}+x_{2}}{2}$, hence $x_{0}, x_{1}\in E$. Repeating this procedure $n$ times, we 
finally infer that $x=x_{n}\in E$.

\Exa{02}{Consider $X=[0,1]$ with the binary operation $\omega:[0,1]^2\to [0,1]$ given by $\omega(x,y):=xy$. Then 
it is easy to check that the $\omega$-extreme sets are of the form: $[p,1], (p,1]$, where $p\in [0,1]$. If $p\in(0,1)$, 
then the $\omega$-extremal hull of $\{p\}$ is $[p,1]$ which does not contain $\omega(p,p)=p^2$, hence $\omega$-extreme
sets may not be $\omega$-convex.}

\Exa{03}{ Let $X=\R$ and let $\omega:\R^2\to \R$ be given by $\omega(x,y):=\min\{x,y\}$. Then the $\omega$-extreme sets 
are of the form: $(p,\infty), [p,\infty)$.} 

\Exa{05}{Let $X=[0,1]$ and let $\omega= \{\omega_{1},\omega_{2}\},$ where $\omega_{1}(x,y):=\frac{x+y}{2}$ and 
$\omega_{2}(x,y):=xy.$ Then the $\omega$-convex subsets of $X$ are of the form $[0,p]$, $[0,p)$, $(0,p)$, $(0,p]$, and 
$\{1\}$, where $p\in [0,1]$. The only $\omega$-extreme sets are $[0,1]$ and $\{1\}.$}

\section{Notions and Properties in Ordered Structures}

In this section we discuss several properties of ordered structures which will be useful in the sequel.
Let $(Y,\leq)$ be a partially ordered set. We start by recalling the following definitions.

An element $u\in Y$ is called the \emph{infimum} (or the \emph{greatest lower bound}) of a nonempty subset $A$ of 
$Y$, written $\inf A$, if
\begin{enumerate}[(a)]
 \item $u$ is a \emph{lower bound} of $A$, i.e., $u\leq y$ holds for all $y\in A$, and
 \item $u$ is the \emph{greatest} lower bound of $A$, i.e., for any lower bound $v$ of $A$, we have $v\leq u$.
\end{enumerate}
The notion of \emph{supremum}, i.e., the \emph{least upper bound} of a nonempty set is defined analogously.

Given another partially ordered set $(Z,\leq)$, a map $\Phi:Y\to Z$ is called an \emph{order preserving map
between $Y$ and $Z$} if, for all $y_1,y_2 \in Y$, the inequality $y_1\leq y_2$ implies $\Phi(y_1)\leq \Phi(y_2)$.
We speak about an \emph{order isomorphism between $Y$ and $Z$} if $\Phi$ is a bijection and, for all 
$y_1,y_2 \in Y$, the condition 
\Eq{*}{
  y_1\leq y_2 \qquad\Longleftrightarrow\qquad \Phi(y_1)\leq \Phi(y_2)
}
holds true. This is equivalent to the property that $\Phi$ and also its inverse $\Phi^{-1}:Z\to Y$ are order 
preserving maps. If $Y=Z$, then $\Phi$ is simply said to be an \emph{order automorphism of $Y$}. The following 
easy-to see lemma shows that the existence of the infimum of a nonempty set is preserved by the action of an order 
isomorphism.

\Lem{CC}{Let $\Phi:Y\to Z$ be an order isomorphism between the partially ordered sets $(Y,\leq)$ and $(Z,\leq)$. Let 
$A\subseteq Y$ be a nonempty lower bounded subset such that $\inf A$ exists. Then $\Phi(A)$ 
is a lower bounded set in $Z$ such that $\inf\Phi(A)$ exists and $\inf\Phi(A)=\Phi(\inf A)$.}

\begin{proof} By the order preserving property of $\Phi$, we obviously have that if $y$ is a lower bound for $A$, 
then $z=\Phi(y)$ is lower bound for $\Phi(A)$. With $y_0:=\inf A$, we get that $z_0:=\Phi(y_0)$ is a lower bound 
for $\Phi(A)$.

Now let $z\in Z$ be any lower bound of $\Phi(A)$. Then, by the order preserving property of $\Phi^{-1}$, $\Phi^{-1}(z)$ 
is a lower bound for $A$, hence $\Phi^{-1}(z)\leq y_0$. This implies that $z\leq\Phi(y_0)$, proving that $z_0$ is the 
largest from among the lower bounds of $\Phi(A)$. Therefore, $z_0$ is the infimum of $\Phi(A)$ and $z_0=\Phi(\inf A)$.
\end{proof}

A set $\L\subseteq Y$ is called a \emph{chain} if any two elements from $\L$ are comparable. 
We say that a partially ordered set $(Y,\leq)$ is \emph{lower chain-complete} if every nonempty lower bounded chain 
has an infimum. In order to describe the most important examples of a lower chain-complete partially ordered set, we 
need to introduce and recall some terminology about partially ordered abelian groups $(Y,+,\leq)$.

A nonempty subsemigroup $S$ of an abelian group $(Y,+)$ is said to be \emph{pointed} and \emph{salient} if $0\in 
S$ and $S\cap (-S)\subseteq\{0\}$, respectively. An arbitrary pointed and salient subsemigroup $S$ of $Y$ induces a 
partial ordering $\leq_{S}$ on $Y$ by letting $x\leq_{S}y$ whenever $y-x \in S$. This partial order is compatible with 
the additive structure of $Y$ in the sense that if $x\leq_{S}y,$ then $x+z\leq_{S} y+z$ for each $z\in Y$. Conversely, 
given a partial order $\leq$ on $Y$ which is compatible with the additive structure of $Y$, the set of nonnegative 
elements of $(Y,+)$, i.e., $S:=\{y\in Y\mid y\geq0\}$ is a pointed and salient subsemigroup of $Y$.

A triple $(Y,+,d)$ is  called a \emph{metric abelian group} if $(Y,+)$ is an abelian group, $(Y,d)$ is a 
metric space and the metric is translation invariant, i.e., $d(x+z,y+z)=d(x,y)$ for all $x,y,z\in Y$. In such a 
case, the metric induces a \emph{pseudo norm} $\|\cdot\|_d:Y\to\R_+$ via the standard definition 
$\|x\|_d:=d(x,0)$. It is easy to see that $\|\cdot\|_d$ is a subadditive and even function. 

If $(Y,+,d)$ is a metric abelian group, then a subsemigroup $S$ is called \emph{additively controllable} if 
there exists a continuous additive function $a:Y\to\R$ such that
\Eq{ac}{
  \|y\|_d \leq a(y)\qquad (y\in S).
}
One can easily see that additively controllable subsemigroups are automatically salient. Indeed, if $y\in 
S\cap(-S)$, then, by \eq{ac}, we get that $\|y\|_d\leq\min(a(y),(a(-y))\leq0$, whence $y=0$ follows.

By the following result, completeness of the metric space $(Y,d)$ and additive controllability of the semigroup of 
nonnegative elements implies the lower chain-completeness of the partially ordered set

\Thm{lcc}{Let $(Y,+,d)$ be a complete metric abelian group and let $S$ be a closed pointed additively controllable 
subsemigroup of $Y$. Then the partially ordered set $(Y,\leq_S)$ is lower chain-complete.}

\begin{proof}
By the controllability assumption, there exist an additive function $a:Y\to\R$ such that \eq{ac} holds.

Let $\Gamma$ be a nonempty set and let $\L:=\{y_{\gamma}\ | \gamma \in \Gamma\}$ be a lower bounded chain in 
$(Y,\leq_S)$ with a lower bound $y_0\in Y$. Since $y_{\gamma}-y_{0}\in S$, therefore we have that $0\leq 
\|y_{\gamma}-y_{0}\|\leq a(y_{\gamma}-y_{0})$ for all $\gamma \in \Gamma.$ This yields that 
\Eq{*}{
\alpha:=\inf_{\gamma \in \Gamma} a(y_{\gamma})\geq a(y_0)>-\infty.
}
By the definition of the infimum, for any $n\in\N$, we can find an element $\gamma_{n}\in \Gamma$ 
such that 
\Eq{*}{
\alpha + \frac{1}{n}>a(y_{\gamma_{n}}).
}
We are now going to show that $(y_{\gamma_{n}})$ is a Cauchy sequence. By the above construction, 
$(a(y_{\gamma_{n}}))$ is a Cauchy sequence (because it converges to $\alpha$). Therefore, for a fixed 
$\varepsilon>0$, there exists $n_{0}\in\N$ such that 
\Eq{*}{
 |a(y_{\gamma_{n}})-a(y_{\gamma_{m}})|<\varepsilon\qquad (n,m\geq n_{0}).
}
Then, in view of the chain property, for $n, m\geq n_{0}$, we have that $y_{\gamma_{n}}-y_{\gamma_{m}}\in 
S\cup(-S)$. Hence, by \eq{ac}, we get
\Eq{*}{
\|y_{\gamma_{n}}-y_{\gamma_{m}}\|_d\leq |a(y_{\gamma_{n}})-a(y_{\gamma_{m}})|<\varepsilon,
}
proving that $(y_{\gamma_{n}})$ is a Cauchy sequence. Let $y_{\star}:=\lim_{n\to\infty}y_{\gamma_{n}}.$ 
It follows from the continuity of $a$ that 
\Eq{*}{
a(y_{\star})=\lim_{n\to \infty}a(y_{\gamma_{n}})=\alpha.
}
We shall show that $y_{\star}=\inf\L.$ First, we prove that $y_{\star}$ is a 
lower bound of the chain $\L.$ Since $y_{\gamma_{n}}-y_{\gamma}\in 
S\cup(-S)$ for all $n\in \N$ and for all $\gamma\in\Gamma$, therefore, by using the closedness of $S$ and 
taking the limit $n\to\infty$, it follows that $y_\star-y_{\gamma}\in S\cup(-S)$ for all $\gamma\in\Gamma$. 
If $y_{\gamma}-y_{\star}\in (-S)\setminus S,$ for some $\gamma \in \Gamma,$ then $y_{\star}-y_{\gamma}\in 
S\setminus \{0\}.$ On account of inequality \eq{ac}, we  obtain
\Eq{*}{
0<\|y_{\star}-y_{\gamma}\|_d\leq a(y_{\star}-y_{\gamma})=a(y_{\star})-a(y_{\gamma})\leq \alpha-\alpha=0,
}
which is a contradiction. Therefore $y_{\gamma}-y_{\star}\in S$, which means that $y_{\star}$ is a lower bound 
of the chain $\L.$ 

If $z\in Y$ is another lower bound of this chain, then $y_{\gamma}-z\in S$ for all $\gamma \in \Gamma.$ In 
particular, $y_{\gamma_{n}}-z\in S,$ for all $n\in \N$. Thus, taking the limit $n\to \infty$, we get 
$y_{\star}-z\in S$. Consequently, $z\leq_{S}y_{\star}$ which means that $y_{\star}=\inf\L.$ The proof of 
theorem is finished.
\end{proof}

With any cone $\K$ in normed space $Y$ we can associate its so-called dual cone ${\K^{\circ}}$, which 
is defined as follows 
\Eq{*}{
 \K^{\circ}:=\{\varphi \in Y^{\star} |\ \varphi(y)\geq 0\ \textrm{for all}\ y\in \K\}.
}
We say that the cone $\K\subseteq Y$ is \emph{sharp} if $\intr({\K^{\circ}})\neq \emptyset.$ 
For the sharp cones, the following useful lemma holds true.

\Lem{lc}{Let $Y$ be a normed space. Then every sharp cone of $Y$ is additively (and therefore linearly) controllable.}

\begin{proof}
Let $\K$ be a sharp cone and let $\varphi \in \intr({\K^{\circ}})$ be a continuous linear functional with 
$\|\varphi\|=1.$ Then there exists a number $r>0$ such that $B(\varphi,r)\subseteq \K^{\circ}.$ We will show that
\Eq{lc}{
  \|y\|\leq \frac1r\varphi(y) \qquad(y\in\K).
}
To prove this, choose an element $y\in \K$ arbitrarily. Then, by a well-known consequence of the Hahn-Banach 
theorem, there exists a linear functional $\psi \in Y^{\star}$ such that
\Eq{*}{
\|y\|=\psi(y) \qquad \textrm{and}\qquad \|\psi\|=1.
}
Then, $\varphi-r\psi\in B(\varphi,r)$, hence
\Eq{*}{
\|y\|=\psi(y)=\frac{1}{r}r\psi(y)=\frac{1}{r}[\varphi(y)-(\varphi-r\psi)(y)]
 \leq \frac{1}{r}\varphi(y),
}
which completes the proof of \eq{lc} showing the linear controllability of $\K$ with the linear functional 
$\frac{1}{r}\varphi$.
\end{proof}

As an immediate consequence of \thm{lcc} and \lem{lc}, we get the following result.

\Cor{01}{Let $(Y,\leq_{\K})$ be a partially ordered vector space, where $Y$ is a Banach space and $\leq_{\K}$ 
is an order generated by a sharp closed convex cone $\K\subseteq Y$. Then $(Y,\leq_{\K})$ is lower chain-complete.}

\begin{proof} Apply \thm{lcc} to the additive group of the vector space $Y$ and to the semigroup $\K$ which,  
 by \lem{lc} is additively controllable.
\end{proof}

We have already seen that sharp cones are always salient. For closed convex cones of finite dimensional 
normed spaces salientness, in fact, is equivalent to sharpness.

\Thm{pc}{Let $Y$ be a finite dimensional normed space. Then every closed convex salient cone of $Y$ is sharp.}

\begin{proof} Assume that $\K$ is a closed convex cone, which is not sharp. Then, $\K^\circ$ is flat, that 
is, it is contained in a proper linear subspace of $Y^\star$. Hence, by the reflexivity of $Y$, there exists 
$y_0\in Y\setminus\{0\}$ such that $\varphi(y_0)=0$ for all $\varphi\in\K^\circ$. On the other hand, in 
view of the so-called bipolar theorem, the convexity and closedness  of $\K$ implies that
\Eq{*}{
  \K=(\K^\circ)^\circ:=\{y\in Y\mid \varphi(y)\geq0\ \textrm{for all}\ \varphi\in \K^\circ\}.
}
Hence $y_0,-y_0\in\K$, which contradicts the salientness of $\K$.
\end{proof}

Another important cone which is sharp is the so-called Lorenz cone. Let $Y$ be a normed space 
and consider the linear space $Y\times \R$ (where as usual, the addition and the scalar multiplication are 
defined coordinatewise). Given a positive number $\varepsilon$, the convex cone $\K_\varepsilon$ defined by 
the formula 
\Eq{*}{
  \K_\varepsilon:=\{(x,t)\in Y\times \R \mid \varepsilon\|x\|\leq t\}
}
is called the \emph{Lorenz cone} (or \emph{ice-cream cone}). 

\Prp{Lc}{Let $Y$ be a normed space. Then, for any positive number $\varepsilon$, the Lorenz cone $\K_\varepsilon$ is a 
sharp closed convex cone in $Y\times \R$.}

\begin{proof} The closedness and convexity of $\K_\varepsilon$ is obvious.
An easy calculation yields that the polar cone of $\K_\varepsilon$ has the form 
\Eq{*}{ 
  \K_\varepsilon^\circ=\{(\varphi,c)\in Y^*\times \R \mid \|\varphi\|+\varepsilon c\leq0\}.
}
Now observe that 
\Eq{*}{
  \intr(\K_\varepsilon^\circ)=\{(\varphi,c)\in Y^*\times \R \mid \|\varphi\|+\varepsilon 
c<0\}\neq \emptyset,
}
which proves that the Lorenz cone is sharp. 
\end{proof}

\section{Convex and affine functions}

In this and in the subsequent sections, we will frequently use the following basic hypothesis which is the 
minimal assumption to formulate our basic definitions and results. 
\begin{enumerate}[(H)]
 \item $X$ is a nonempty set and $(Y,\leq)$ is a partially ordered set,
  $\Gamma$ is a nonempty set, $n:\Gamma\to\N$ is an arity function and
       $ \omega=\{\omega_\gamma:X^{n(\gamma)}\to X\mid\gamma\in\Gamma\}$ and  
       $\Omega=\{\Omega_\gamma:Y^{n(\gamma)}\to Y\mid\gamma\in\Gamma\}$ are two given families of operations. 
\end{enumerate}

A family of operations $\omega=\{\omega_{\gamma} \mid \gamma\in \Gamma \}$ is said to be  \emph{a pairwise mutually 
distributive} if for all $\gamma, \beta \in \Gamma,\ k\in \{1,2,\dots,n(\gamma)\}$ and all 
$x_1,\dots,x_{k-1}$, $x_{k+1},\dots,x_{n(\gamma)}$, $y_1,\dots,y_{n(\beta)}\in X$ 
\Eq{*}{
\omega_{\gamma}(x_1,\dots,x_{k-1},&\omega_{\beta}(y_1,\dots,y_{n(\beta)}),x_{k+1},\dots,x_{n(\gamma)}) \\ 
&=\omega_{\beta}(\omega_{\gamma}(x_1,\dots,x_{k-1},y_1,x_{k+1},\dots,x_{n(\gamma)}),\dots, \\
& \qquad\qquad\omega_{\gamma}(x_1,\dots,x_{k-1},y_{n(\beta)},x_{k+1},\dots,x_{n(\gamma)})).
}
We say that a family of operations $\omega=\{\omega_{\gamma} \mid \gamma\in \Gamma \}$ is \emph{reflexive} if, for all 
$\gamma \in \Gamma$,  
\Eq{*}{
 \omega_{\gamma}(x,\dots,x)=x,\quad x\in X.
}
 
Under the hypothesis (H), given an $\omega$-convex set $D\subseteq X$, we say that $f:D\to Y$ is 
\emph{$(\omega,\Omega)$-convex on $D$} if it satisfies the functional inequality
\Eq{*}{
 f\big(\omega_\gamma(x_1,\dots,x_{n(\gamma)})\big)
    \leq \Omega_\gamma\big(f(x_1),\dots,f(x_{n(\gamma)})\big) \quad 
    (\gamma\in\Gamma,\,x_1,\dots,x_{n(\gamma)}\in D).
}
If $f$ satisfies the reversed inequality 
\Eq{*}{
 \Omega_\gamma\big(f(x_1),\dots,f(x_{n(\gamma)})\big) 
   \leq  f\big(\omega_\gamma(x_1,\dots,x_{n(\gamma)})\big)
   \quad (\gamma\in\Gamma,\,x_1,\dots,x_{n(\gamma)}\in D),
} 
then we say that it is \emph{$(\omega,\Omega)$-concave on $D$}.
Finally, a function $f$ is called \emph{$(\omega,\Omega)$-affine on $D$} if it satisfies the functional equation
\Eq{*}{
 f\big(\omega_\gamma(x_1,\dots,x_{n(\gamma)})\big)
    = \Omega_\gamma\big(f(x_1),\dots,f(x_{n(\gamma)})\big) \quad 
(\gamma\in\Gamma,\,x_1,\dots,x_{n(\gamma)}\in D).
}
Trivially, a function is $(\omega,\Omega)$-affine if and only if it is $(\omega,\Omega)$-convex and 
$(\omega,\Omega)$-concave.

The basic properties of $(\omega,\Omega)$-convexity with respect to the pointwise supremum and infimum are established 
in the following results.

\Thm{sup}{Assume that the hypothesis (H) holds and, for all $\gamma\in\Gamma$, 
the operation $\Omega_\gamma$ is nondecreasing with respect to each of its variables. Let $D\subseteq X$ be an 
$\omega$-convex set, $\Delta$ be a nonempty set, $\mathscr{F}=\{f_\delta:D\to Y\mid\delta\in\Delta\}$ be a family of 
$(\omega,\Omega)$-convex functions on $D$ and $f:D\to Y$.
\begin{enumerate}
\item If either $f$ satisfies
\Eq{fsup}{
  f(x):=\sup\{f_\delta(x)\mid\delta\in\Delta\} \qquad(x\in D),
}
\item or $\mathscr{F}$ is a chain with respect to the pointwise ordering, for all $\gamma\in\Gamma$, 
the operation $\Omega_\gamma$ is an order isomorphism with respect to each of its variables, and $f$ satisfies
\Eq{finf}{
  f(x):=\inf\{f_\delta(x)\mid\delta\in\Delta\} \qquad(x\in D),
}
\end{enumerate}
then $f$ is $(\omega,\Omega)$-convex on $D$.}

\begin{proof} First assume that $f$ is given by \eq{fsup}. To prove its $(\omega,\Omega)$-convexity,
let $\gamma\in\Gamma$ and $x_1,\dots,x_{n(\gamma)}\in D$ be arbitrary. Then, by the $(\omega,\Omega)$-convexity of 
$f_\delta$ and by the monotonicity property of $\Omega_\gamma$, for all $\delta\in\Delta$, we get
\Eq{*}{
  f_\delta\big(\omega_\gamma(x_1,\dots,x_{n(\gamma)})\big)
    \leq \Omega_\gamma\big(f_\delta(x_1),\dots,f_\delta(x_{n(\gamma)})\big)
    \leq \Omega_\gamma\big(f(x_1),\dots,f(x_{n(\gamma)})\big).
}
Upon taking the supremum of the left hand side of this inequality with respect to $\delta\in\Delta$, it follows that
\Eq{*}{
  f\big(\omega_\gamma(x_1,\dots,x_{n(\gamma)})\big)
    =\sup_{\delta\in\Delta}f_\delta\big(\omega_\gamma(x_1,\dots,x_{n(\gamma)})\big)
    \leq \Omega_\gamma\big(f(x_1),\dots,f(x_{n(\gamma)})\big),
}
which shows that $f$ is $(\omega,\Omega)$-convex.

Secondly, assume that $\mathscr{F}$ is a chain and $f$ satisfies \eq{finf}. To verify the $(\omega,\Omega)$-convexity 
of $f$, let $\gamma\in\Gamma$ and $x_1,\dots,x_{n(\gamma)}\in D$ be fixed and let 
$\delta_1,\dots,\delta_{n(\gamma)}\in\Delta$ be arbitrary. Using that $\mathscr{F}$ is a chain, the existence of 
an index $\delta_*\in\{\delta_1,\dots,\delta_{n(\gamma})\}$ can be established such that $f_{\delta_*}(x)\leq 
f_{\delta_i}(x)$ holds for all $x\in D$ and $i\in\{1,\dots,n(\gamma)\}$. Then, by the $(\omega,\Omega)$-convexity of 
$f_{\delta_*}$ and by the monotonicity property of the operation $\Omega_\gamma$, we get
\Eq{*}{
  f\big(\omega_\gamma(x_1,\dots,x_{n(\gamma)})\big)
    &\leq f_{\delta_*}\big(\omega_\gamma(x_1,\dots,x_{n(\gamma)})\big) \\
    &\leq \Omega_\gamma\big(f_{\delta_*}(x_1),\dots,f_{\delta_*}(x_{n(\gamma)})\big) \\
    &\leq \Omega_\gamma\big(f_{\delta_{1}}(x_1),\dots,f_{\delta_{n(\gamma)}}(x_{n(\gamma)})\big)
}
for all $\delta_1,\dots,\delta_{n(\gamma)}\in\Delta$. Using that $\Phi_\gamma$ is an order isomorphism in its first 
variable, for all $\delta_2,\dots,\delta_{n(\gamma)}\in\Delta$, we get
\Eq{*}{
  f\big(\omega_\gamma(x_1,\dots,x_{n(\gamma)})\big)
    &\leq\inf_{\delta_1\in\Gamma}\Omega_\gamma\big(f_{\delta_1}(x_1),
         f_{\delta_2}(x_2),\dots,f_{\delta_{n(\gamma)}}(x_{n(\gamma)})\big)\\
    &\leq \Omega_\gamma\Big(\inf_{\delta_1\in\Gamma}f_{\delta_{1}}(x_1),
         f_{\delta_2}(x_2),\dots,f_{\delta_{n(\gamma)}}(x_{n(\gamma)})\Big)\\
  &=\Omega_\gamma\Big(f(x_1),
         f_{\delta_2}(x_2),\dots,f_{\delta_{n(\gamma)}}(x_{n(\gamma)})\Big).
}
(In the case when $n(\gamma)=1$, the above inequalities can easily be adjusted.) Repeating this step and taking the 
infimum with for $\delta_2,\dots,\delta_{n(\gamma)}$, respectively, we shall arrive at the inequality  
\Eq{*}{
  f\big(\omega_\gamma(x_1,\dots,x_{n(\gamma)})\big)
   \leq\Omega_\gamma\Big(f(x_1),f(x_2),\dots,f(x_{n(\gamma)})\Big),
}
which proves that $f$ is $(\omega_\gamma,\Omega_\gamma)$-convex. This completes the proof of the 
$(\omega,\Omega)$-convexity of $f$.
\end{proof}

\Cor{sup}{In addition to the assumption (H), suppose that, for all $\gamma\in\Gamma$, the operation 
$\Omega_\gamma$ is nondecreasing with respect to each of its variables. Let $D\subseteq X$ be an 
$\omega$-convex set, let $\Delta$ be a nonempty set, let $\{g_\delta:D\to Y\mid\delta\in\Delta\}$ be 
a family of $(\omega,\Omega)$-affine functions on $D$ and assume that $f:D\to Y$ satisfies 
\Eq{sup}{
  f(x)=\sup\{g_\delta(x)\mid\delta\in\Delta\} \qquad(x\in D).
}
Then $f$ is $(\omega,\Omega)$-convex on $D$.}

\begin{proof} Since $(\omega,\Omega)$-affine functions are automatically $(\omega,\Omega)$-convex, therefore
the first part of \thm{sup} yields the statement.
\end{proof}

Our first main result establishes the affine extension of a function which is dominated by a convex one.

\Thm{MT1}{In addition to hypothesis (H) above, assume that
\begin{enumerate}[(H1)]
 \item $(Y,\leq)$ is a lower chain-complete partially ordered set.
 \item The family $\omega$ consists of pairwise mutually distributive operations.
 \item The family $\Omega$ consists of pairwise mutually distributive operations 
such that, for all $\gamma\in\Gamma$, the operation $\Omega_\gamma$ is an order automorphism in each of its 
variables. 
\end{enumerate}
Let $f:X\to Y$ be an $(\omega,\Omega)$-convex function and let $D\subseteq X$ be a nonempty $\omega$-convex subset of 
$X$ such that $\ext_\omega(D)=X$ and $f|_D$ is $(\omega,\Omega)$-affine on $D$. Then there exists an 
$(\omega,\Omega)$-affine function $g:X\to Y$ such that $g\leq f$ and $g|_D=f|_D$.}

\begin{proof} For the proof of the theorem, consider the following collection of functions mapping $X$ into $Y$:
\Eq{*}{
 \G:=\{g:X\to Y\mid \mbox{$g$ is $(\omega,\Omega)$-convex, $g\leq f$ and } g|_D=f|_D\}.
}
Our aim is to verify that $\G$ contains an $(\omega,\Omega)$-affine element.

First observe that $\G$ is not empty because $f\in\G$ trivially holds. Observe that the family $\G$ can be partially 
ordered using the partial order of $Y$ by letting $g\leq h$ if and only if $g(x)\leq h(x)$ for all $x\in X$.
By Zorn's Lemma, there exists a maximal chain $\{g_\delta\in\G\mid\delta\in\Delta\}$ in the partially ordered set 
$(\G,\leq)$. We are going to prove that the infimum of this chain exists and is an $(\omega,\Omega)$-affine 
function.

Denote by $E\subseteq X$ the set of those points $x$ such that $\{g(x)\mid g\in\G\}$ is lower bounded. Because, for 
$g\in\G$, we have that $g|_D=f|_D$, hence $D\subseteq E$. We show that $E$ is $\omega$ extreme. To see this, let 
$\gamma\in\Gamma$ and let $(x_1,\dots,x_{n(\gamma)})\in\omega_\gamma^{-1}(E)$. This means that 
$\omega_\gamma(x_1,\dots,x_{n(\gamma)})$ is in $E$. Let $y_0$ denote a lower bound for the set 
$\{g\big(\omega_\gamma(x_1,\dots,x_{n(\gamma)})\big)\mid g\in\G\}$. Then, for $g\in\G$, by the 
$(\omega,\Omega)$-convexity of $g$ and by the inequality $g\leq f$, we get that
\Eq{*}{
  y_0\leq g\big(\omega_\gamma(x_1,\dots,x_{n(\gamma)})\big)
  &\leq \Omega_\gamma\big(g(x_1),\dots,g(x_{n(\gamma)})\big)\\
  &\leq \Omega_\gamma\big(f(x_1),\dots,g(x_i),\dots,f(x_{n(\gamma)})\big)
}
if $i\in\{1,\dots,n(\gamma)\}$. In view of the order automorphism property of $\Omega_\gamma$ in its $i$th variable, it 
follows that the set $\{g(x_i)\mid g\in\G\}$ is lower bounded, i.e., $x_i\in E$ for all $i\in\{1,\dots,n(\gamma)\}$.
This proves that $E$ is $\omega$-extreme, whence the assumption $\ext_\omega(D)=X$ and the inclusion $D\subseteq E$ 
imply that $E=X$.

Therefore, for all $x\in X$, the chain $\{g_\delta(x)\in\G\mid\delta\in\Delta\}$ is lower bounded. Applying the lower 
chain completeness of $Y$, it follows that the set $\{g_\delta(x)\in\G\mid\delta\in\Delta\}$ has an infimum, which we 
will denote by $g_0(x)$. The function $g_0:X\to Y$ so defined is $(\omega,\Omega)$-convex by the second 
assertion of \thm{sup}. 

To complete the proof, it is enough to show that $g_0$ is an $(\omega,\Omega)$-affine on $X$. Because, for all 
$\delta\in\Delta$, $g_\delta$ equals $g$ on $D$, therefore $g_0$ is also equal to $g$ on $D$, and hence it is 
$(\omega,\Omega)$-affine on $D$.

Now, fix $\gamma \in \Gamma$ arbitrarily. We prove by induction on $k\in \{0,\dots,n(\gamma)\}$ that the equality
\Eq{ind}{
g_0\big(\omega_{\gamma}(x_1,\dots, x_k,&y_{k+1},\dots, y_{n(\gamma)})\big)\\
&=\Omega_{\gamma}\big(g_0(x_1),\dots,g_0(x_k),g_0(y_{k+1}),\dots,g_0(y_{n(\gamma)})\big)
}
holds for all $x_1,\dots, x_k\in X$ and $y_{k+1},\dots, y_{n(\gamma)}\in D$. (We accept here the convention 
that if $k=0$ (resp.\ $k=n(\gamma)$) then the $x_i$s (resp.\ $y_j$s) are missing.) The 
statement is obvious for $k=0$ due to the $(\omega,\Omega)$-affine of $g_0$ on $D$.

Assume \eq{ind} for some $k\in \{0,\dots,n(\gamma)-1\}$ and for all $x_1,\dots, x_k\in X,\ y_{k+1},\dots, 
y_{n(\gamma)}\in D$. Fix $x_1,\dots,x_k\in X$ and $ y_{k+2},\dots,y_{n(\gamma)}\in D$ arbitrarily. 
Because $\Omega_{\gamma}$ is an order automorphism with respect to its $(k+1)$st variable, thus there exists a 
uniquely determined function $u:X\to Y$ such that, for all $x_{k+1}\in X$, we have 
\Eq{fff}{
  \Omega_\gamma\big(g_0(x_1),\dots,g_0(x_k),&u(x_{k+1}),g_0(y_{k+2}),\dots,g_0(y_{n(\gamma)})\big)\\
  &= g_0\big(\omega_\gamma(x_1,\dots,x_k,x_{k+1},y_{k+2},\dots,y_{n(\gamma)})\big).
}
By using the $(\omega,\Omega)$-convexity of $g_0$, we infer that
\Eq{*}{
\Omega_\gamma\big(g_0(x_1),\dots,&g_0(x_k),u(x_{k+1}),g_0(y_{k+2}),\dots,g_0(y_{n(\gamma)})\big) \\
&\leq \Omega_\gamma\big(g_0(x_1),\dots,g_0(x_k),g_0(x_{k+1}),g_0(y_{k+2}),\dots,g_0(y_{n(\gamma)})\big).
}
The order automorphism property of $\Omega_{\gamma}$ with respect to its $(k+1)$st variable implies that $u\leq 
g_0$ on $X$.

We will show that $u\in \G$. First, observe that $u$ is an $(\omega,\Omega)$-convex map. Indeed, let $\beta\in \Gamma$ 
and 
$z_1,\dots,z_{n(\beta)}\in X$. Then, using \eq{fff} for $x_{k+1}:=\omega_{\beta}(z_1,\dots,z_{n(\beta)})$, then 
assumption (H2), next the $(\omega,\Omega)$-convexity of $g_0$, then \eq{fff} again, finally the assumption (H3), we get
\Eq{*}{
&\Omega_{\gamma}\big(g_0(x_1),\dots,g_0(x_k),u(\omega_{\beta}(z_1,\dots,z_{n(\beta)})),
       g_0(y_{k+2}),\dots,g_0(y_{n(\gamma)})\big) \\ 
&=g_0\big(\omega_\gamma(x_1,\dots,x_k,\omega_\beta(z_1,\dots,z_{n(\beta)}),y_{k+2},\dots,y_{n(\gamma)})\big) \\ 
&=g_0\big(\omega_\beta(\omega_\gamma(x_1,\dots,x_k,z_1,y_{k+2},\dots,y_{n(\gamma)}),\dots,
\omega_\gamma(x_1,\dots,x_k,z_{n(\beta)},y_{k+2},\dots,y_{n(\gamma)}))\big) \\ 
&\leq\Omega_{\beta}\big(g_0(\omega_\gamma(x_1,\dots,x_k,z_1,y_{k+2},\dots,y_{n(\gamma)})),\dots, 
g_0(\omega_\gamma(x_1,\dots,x_k,z_{n(\beta)},y_{k+2},\dots,y_{n(\gamma)}))\big) \\ 
&=\Omega_\beta\big(\Omega_\gamma(g_0(x_1),\dots,u(z_1),\dots g_0(y_{n(\gamma)})),\dots, 
\Omega_\gamma(g_0(x_1),\dots,u(z_{n(\beta)}),\dots g_0(y_{n(\gamma)}))\big) \\ 
&=\Omega_\gamma\big(g_0(x_1),\dots,g_0(x_k),\Omega_\beta(u(z_1),\dots,u(z_{n(\beta)})),
    g_0(y_{k+2}),\dots,g_0(y_{n(\gamma)})\big).
}
Using again the order automorphism property of $\Omega_{\gamma}$ with respect to its $(k+1)$st variable, we obtain 
that 
\Eq{*}{
u(\omega_\beta(z_1,\dots,z_{n(\beta)}))\leq \Omega_\beta(u(z_1),\dots,u(z_{n(\beta)})),
} 
which completes the proof of the $(\omega,\Omega)$-convexity of $u$. 

Now, let us observe that $u_{|D}=f_{|D}$. Indeed, using the inductive assumption, that is the validity of \eq{ind} for 
$k$, and also formula \eq{fff}, for all $y\in D$, we obtain 
\Eq{*}{
\Omega_\gamma\big(g_0(x_1),&\dots,g_0(x_k),u(y),g_0(y_{k+2}),\dots,g_0(y_{n(\gamma)})\big) \\
&=g_0(\omega_\gamma(x_1,\dots,x_k,y,y_{k+2},\dots,y_{n(\gamma)})\\
&= \Omega_\gamma\big(g_0(x_1),\dots,g_0(x_k),g_0(y),g_0(y_{k+2}),\dots,g_0(y_{n(\gamma)})\big).
} 
Therefore, the order automorphism property of $\Omega_{\gamma}$ with respect to its $(k+1)$st variable, 
yields that $u(y)=g_0(y)$ for all $y\in D$. We have shown that $u\in \G$. On the other hand $u\leq g_0$ 
then, in view of the minimality of $g_0$, it follows that $u=g_0$. 
Hence 
\Eq{*}{
g_0(\omega_\gamma(x_1,\dots,&x_k,x_{k+1},y_{k+2},\dots,y_{n(\gamma)}))\\
&=\Omega_\gamma(g_0(x_1),\dots,g_0(x_k),g_0(x_{k+1}),g_0(y_{k+2}),\dots,g_0(y_{n(\gamma)})),
}
for all $x_1,\dots,x_{k+1}\in X,\ y_{k+2},\dots,y_{n(\gamma)}\in D$ which finishes the proof of \eq{ind} for all 
$k\in\{0,\dots,n(\gamma)\}$. 

Finally, applying \eq{ind} for $k=n(\gamma)$, we obtain that $g_0$ is $(\omega_\gamma,\Omega_\gamma)$-affine. 
Since $\gamma\in\Gamma$ was arbitrary, this yields that $g_0$ is $(\omega,\Omega)$-affine, which was to be proved.
\end{proof}

The following consequence of the above theorem is a support theorem which, in some sense, reverses the statement of 
\cor{sup}. Here, we have to assume that the operations involved are reflexive.

\Cor{MT2}{In addition to hypothesis (H) above, assume that
\begin{enumerate}[(H1+)]
 \item $(Y,\leq)$ is a lower chain-complete partially ordered set.
 \item The family $\omega$ consists of reflexive and pairwise mutually distributive operations.
 \item The family $\Omega$ consists of reflexive and pairwise mutually distributive operations 
such that, for all $\gamma\in\Gamma$, the operation $\Omega_\gamma$ is an order automorphism in each of its 
variables. 
\end{enumerate}
Let $f:X\to Y$ be an $(\omega,\Omega)$-convex function. Then, for all $\omega$-interior point 
$p\in X$, there exists an $(\omega,\Omega)$-affine function $g:X\to Y$ such that $g\leq f$ and $g(p)=f(p)$.}

\begin{proof}
Put $D:=\{p\}$. Obviously, due to the reflexivity property of each $\omega_{\gamma} \in \omega$ the set $D$ is 
$\omega$-convex. The reflexivity of the operations $\omega_\gamma$ and $\Omega_\gamma$ implies that $f|_D$ is 
$(\omega,\Omega)$-affine. Now, to finish the proof, it is enough to use the \thm{MT1}.
\end{proof}

In the subsequent result we apply \thm{MT1} and \cor{MT2} to various situations when the operations are given in terms 
of additive maps. 

\Cor{MTA}{Let $(X,+)$ be an abelian semigroup, and let $(Y,+,d)$ be a complete metric abelian group equipped with an 
ordering $\leq_S$ generated by a closed pointed additively controllable semigroup $S\subseteq Y$.  
Let $f:X\to Y$ be subadditive, i.e., assume that, for all $x,y\in X$, 
\Eq{MTA}{
f(x+y)\leq_{S} f(x)+f(y)
}
holds. Assume that $p\in X$ possesses the following two properties: 
\begin{enumerate}[(i)]
 \item for all $n\in\N$, $f(np)=nf(p)$;
 \item for all $x\in X$, there exist $y\in X$ and $n\in\N$ such that $x+y=np$.
\end{enumerate}
Then there exists an additive function $g:X\to Y$ such that $g\leq_S f$ and $g(p)=f(p)$.}

\begin{proof}
Let $\Gamma=\{1\}$, $n(1)=2$, $\omega=\{\omega_1\}$ and $\Omega=\{\Omega_1\}$, where the operations $\omega_1:X^2\to X$ 
and $\Omega_1:Y^2\to Y$ are given by the formulas: 
\Eq{*}{
  \omega_1(x_1,x_2):=x_1+x_2,\qquad  
  \Omega_1(y_1,y_2):=y_1+y_2. 
}
These operations are obviously autodistributive. Furthermore, a function $f:X\to Y$ is $(\omega,\Omega)$-convex 
(resp.\ $(\omega,\Omega)$-affine) if and only if $f$ is subadditive (resp.\ additive). 

Define the set $D\subseteq X$ by $D:=\{np\mid n\in\N\}$. Then, $D$ is closed under addition, therefore it is 
$\omega$-convex. By assumption (i), $f$ is additive on $D$, which implies that $f$ is $(\omega,\Omega)$-affine on $D$. 
Let $E\subseteq X$ be any $\omega$-extreme set containing $D$. By property (ii), for every $x\in X$, there exists 
$y\in X$ such that $\omega_1(x,y)\in D\subseteq E$. Thus, the $\omega$-extremality of $E$ implies that $(x,y)\in E^2$, 
whence $x\in E$ follows. Therefore, we get that $E=X$ proving that $\ext_\omega(D)=X$. 

In view of \thm{MT1}, there exists an $(\omega,\Omega)$-affine, (i.e., additive) function $g:X\to Y$ such that 
$g\leq_S f$ and $g(p)=f(p)$. The proof is complete.
\end{proof}

\Cor{MTB}{Let $X$ be a convex cone of a linear space, let $Y$ be a Banach space equipped with an ordering $\leq_K$ 
generated by a sharp closed cone $\K\subseteq Y$. Let $f:X\to Y$ be sublinear, i.e., assume that, for all $x,y\in X$ 
and 
$t,s>0$,
\Eq{MTB}{
f(tx+sy)\leq_{S} tf(x)+sf(y)
}
holds. Assume that $p\in X$ possesses the following two properties: 
\begin{enumerate}[(i)]
 \item for all $t>0$, $f(tp)=tf(p)$;
 \item for all $x\in X$, there exist $y\in X$ and $t>0$ such that $x+y=tp$.
\end{enumerate}
Then there exists an additive and positively homogeneous function $g:X\to Y$ such that $g\leq_S f$ and $g(p)=f(p)$.}

\begin{proof}
Let $\Gamma=\{(t,s)\mid t,s>0\}$, $n(t,s)=2$, $\omega=\{\omega_{(t,s)}\mid t,s>0\}$ and $\Omega=\{\Omega_{(t,s)}\mid 
t,s>0\}$, where the operations $\omega_{(t,s)}:X^2\to X$ and $\Omega_{(t,s)}:Y^2\to Y$ are given by the formulas: 
\Eq{*}{
  \omega_{(t,s)}(x_1,x_2):=tx_1+sx_2,\qquad  
  \Omega_{(t,s)}(y_1,y_2):=ty_1+sy_2. 
}
It is easy to check that these operations are distributive with respect to each other. Furthermore, a function $f:X\to 
Y$ is $(\omega,\Omega)$-convex (resp.\ $(\omega,\Omega)$-affine) if and only if $f$ is subadditive (resp.\ additive) 
and positively homogeneous. 

Define the set $D\subseteq X$ by $D:=\{tp\mid t>0\}$. Then, $D$ is closed under addition and multiplication by 
positive scalars, therefore it is $\omega$-convex. By assumption (i), $f$ is additive and positively homogeneous on 
$D$, 
which implies that $f$ is $(\omega,\Omega)$-affine on $D$. Let $E\subseteq X$ be any $\omega$-extreme set containing 
$D$. By property (ii), for every $x\in X$, there exists $y\in X$ such that $\omega_{(1,1)}(x,y)\in D\subseteq E$. Thus, 
the $\omega$-extremality of $E$ implies that $(x,y)\in E^2$, whence $x\in E$ follows. Therefore, we get that $E=X$ 
proving that $\ext_\omega(D)=X$. 

In view of \thm{MT1}, there exists an $(\omega,\Omega)$-affine, (i.e., additive and positively homogeneous) function 
$g:X\to Y$ such that $g\leq_S f$ and $g(p)=f(p)$. The proof is complete.
\end{proof}

For the formulation of the conditions of the subsequent result, we first recall some well-known concepts. 
An abelian group $(G,+)$ is called \textit{uniquely 2-divisible} if, for every $x\in G$, there exists a unique element 
$y\in G$ such that $2y=x$. This element $y$ will be denoted $\frac12x$. The expression $\frac1{2^n}x$ is defined by 
induction with respect to $n\in\N$. Let $X$ be a subset of uniquely 2-divisible abelian group $(G,+)$. $X$ is said to 
be \textit{midconvex} if, for all $x,y\in X$, the midpoint $\frac{1}{2}(x+y)$ also belongs to $X$ (cf.\ 
\cite{JarPal15}). It easily follows by induction, that if $X$ is midconvex then, it is closed under diadic rational 
convex combinations, that is, for all $x,y\in X$ and for all $n\in\N$, $k\in\{0,1,\dots,2^n\}$, the element 
$\frac{k}{2^n}x+(1-\frac{k}{2^n})y$ is contained in $X$. We say that $p$ is a \textit{relative algebraic interior 
point} 
of the set $X$ if, for all $x\in X$, there exists $n\in\N$ such that $p+\frac{1}{2^n}(p-x)\in X$. The set of relative 
algebraic interior points of $X$ will be denoted by $\ri(X)$.

\Lem{ri}{Let $X$ be a midconvex subset of a uniquely 2-divisible abelian group $(G,+)$. Let $a:G\to G$ be an additive 
map and define the operation $\omega:G^2\to G$ by $\omega(x,y):=a(x)+y-a(y)$. Assume that $X$ is $\omega$-convex, 
i.e., $\omega(X^2)\subseteq X$. Then 
\Eq{*}{
   \ri(X)\subseteq\intr_\omega(X).
}}

\begin{proof} Let $p\in\ri(X)$ be arbitrarily fixed. Denote the $\omega$-extreme hull of $\{p\}$ by $E$. In order to 
prove that $p\in\intr_\omega(X)$, we have to show that $E=X$. Let $x\in X$ be arbitrary. By $p\in\ri(X)$, there exists 
$n\in\N$ such that $p+\frac1{2^n}(p-x)\in X$. Define the sequence $x_{-2},x_{1},x_0,\dots,x_{2^{n+1}}$ as follows:
\Eq{ri1}{
  x_{2k}&:=\tfrac{k}{2^n}x+(1-\tfrac{k}{2^n})p \qquad &&(k\in\{-1,0,\dots,2^n\}),\\
  x_{2k-1}&:=\omega(x_{2k-2},x_{2k}) \qquad &&(k\in\{0,\dots,2^n\}).
}
Obviously, $x_0=p$ and $x_{2^{n+1}}=x$. 
Due to $p+\frac1{2^n}(p-x)\in X$, we have that $x_{-2}\in X$. The midconvexity of $X$ implies that $x_{2k}\in X$ for 
all $k\in\{0,\dots,2^n\}$. On the other hand, by the $\omega$-convexity of $X$, it follows that $x_{2k-1}\in X$ for 
all $k\in\{0,\dots,2^n\}$. Therefore all members of the sequence $x_{-2},x_{1},x_0,\dots,x_{2^{n+1}}$ belong to $X$.
We are now going to show that
\Eq{ri2}{
  x_{2k}:=\omega(x_{2k+1},x_{2k-1}) \qquad (k\in\{0,\dots,2^n-1\}).
}
For brevity, denote the additive mapping $\id_G-a$ by $b$. Then, the operation $\omega$ is given 
by $\omega(x,y)=a(x)+b(y)$ and we also have the following two easy-to-see properties of $b$: 
\Eq{ab}{
  a+b=\id_G\qquad \mbox{and}\qquad a\circ b=b\circ a.
}
Denote the element $\frac1{2^n}(x-p)$ by $u$. Then, for $k\in\{0,\dots,2^n-1\}$, we have
\Eq{u}{
  x_{2k\pm 2}=\tfrac{k\pm1}{2^n}x+(1-\tfrac{k\pm1}{2^n})p
  =\tfrac{k}{2^n}x+(1-\tfrac{k}{2^n})\pm\tfrac1{2^n}(x-p)
  =x_{2k}\pm u.
}
Therefore, using \eq{ri1}, \eq{u} and finally the identities of \eq{ab}, we get
\Eq{*}{
 \omega(x_{2k+1},x_{2k-1})
  &=a(x_{2k+1})+b(x_{2k-1}) 
  =a\big(\omega(x_{2k},x_{2k+2})\big)+b\big(\omega(x_{2k-2},x_{2k})\big)\\
  &=a\big(a(x_{2k})+b(x_{2k+2})\big)+b\big(a(x_{2k-2})+b(x_{2k}))\big)\\
  &=a\big(a(x_{2k})+b(x_{2k}+u)\big)+b\big(a(x_{2k}-u)+b(x_{2k}))\big)\\
  &=(a\circ a+a\circ b+b\circ a+b\circ b)(x_{2k})+(a\circ b-b\circ a)(u)\\
  &=\big((a+b)\circ (a+b)\big)(x_{2k})=x_{2k}.
}

Using that $\omega(x_1,x_{-1})=x_0=p\in E$, it follows that $x_1\in E$. Next, applying that $\omega(x_0,x_2)=x_1\in E$, 
we obtain that $x_2\in E$. Using the second equality in \eq{ri1} and equation \eq{ri2} alternately, we infer that 
$x_k$ is in $E$ for all $k\in\{0,\dots,2^{n+1}\}$. In particular, $x$ is contained in $E$, which completes the proof of 
the inclusion $X\subseteq E$.
\end{proof}

\Thm{MT2}{Let $X$ be a midconvex subset of a uniquely 2-divisible abelian group $(G,+)$, and let $(Y,+,d)$ be a 
complete metric abelian group equipped with an ordering $\leq_S$ generated by a closed pointed additively controllable 
semigroup $S\subseteq Y.$ Moreover, assume that $n\geq2$ and $a_{1},\ldots, a_{n}:G\to G$ and $A_{1},\ldots,A_{n}:Y\to 
Y$  are two families of additive maps with the following additional properties: 
\begin{enumerate}[(i)]
\item $a_i\circ a_j=a_j\circ a_i$ and $A_i\circ A_j=A_j\circ A_i,$\ for all $i, j =1,\dots,n$; 
\item $a_{1}+\dots+a_{n}=\id_G$ and $A_{1}+\dots+A_{n}=\id_Y$; 
\item $a_{1}(X)+\ldots+a_{n}(X)\subseteq X$; 
\item $A_i$ is bijective with $A_i(S)=S$ for all $i\in\{1,\dots,n\}.$ 
\end{enumerate}
Let $f:X\to Y$ satisfy, for all $x_1,\ldots,x_n\in X$, the following convexity type inequality
\Eq{MT2c}{
f\big(a_{1}(x_{1})+\dots +a_{n}(x_{n})\big)\leq_{S} A_{1}\big(f(x_{1})\big) + \dots +A_{n}\big(f(x_{n})\big).
}
Then, for every $p\in \ri(X)$, there exists a function $g:G\to Y$ such that $g\leq_S f$, $g(p)=f(p)$ and, for all
$x_1,\ldots,x_n\in X$, the following functional equation holds:
\Eq{MT2a}{
g\big(a_{1}(x_{1})+\dots +a_{n}(x_{n})\big)= A_{1}\big(g(x_{1})\big) + \dots +A_{n}\big(g(x_{n})\big).
}}

\begin{proof}
First, observe that on the account of \cor{01}, the space $(Y,\leq_S)$ is a lower chain complete partially ordered set. 
Let $\Gamma=\{1\},\ n(1)=n,\ \omega=\{\omega_1\}$ and $\Omega=\{\Omega_1\}$, where the operations $\omega_1:G^n\to G$ 
and $\Omega_1:Y^n\to Y$ are given by the formulas: 
\Eq{*}{
  \omega_1(x_1,\dots,x_n)&:=a_1(x_1)+\cdots +a_n(x_n),\\\ 
  \Omega_1(y_1,\dots,y_n)&:=A_1(y_1)+\cdots +A_n(y_n). 
}
These operations are autodistributive due to the pairwise commutativity property of the families 
$\{a_1,\dots,a_n\}$ and $\{A_1,\dots,A_n\}$ postulated in (i). The reflexivity of both operations follows 
from the assumption (ii). In view of property (iii), we have that $\omega_1(X^n)\subseteq X$, that is, $X$ is 
$\omega$-convex.

The operation $\Omega_1$ is an order automorphism in each of its variables, since  the additive maps $A_1,\dots,A_n$ 
are bijective with condition $A_i(S)=S$ for all $i\in \{1,\dots,n\}$. To see this, let $i\in \{1,\dots,n\}$ and 
$y_1,\dots,y_{i-1},y_{i+1},\dots,y_n\in Y$ be fixed. The map $A_i$ being a bijection of $Y$ onto itself, it follows that
\Eq{*}{
  y\mapsto \Omega_1(y_1,\dots,y_{i-1},y,y_{i+1},\dots,y_n)=A_i(y)+\sum_{j\in\{1,\dots,n\}\setminus\{i\}}A_j(y_j)
}
is also a bijection of $Y$ onto itself. On the other hand, applying condition $A_i(S)=S$, for all $y',y''\in X$,
\Eq{*}{
  y'&\leq_S y'' \\ 
  &\Leftrightarrow\, y''-y'\in S \\
  &\Leftrightarrow\, A_i(y'')-A_i(y')=A_i(y''-y')\in A_i(S)=S \\
  &\Leftrightarrow\, A_i(y')\leq_S A_i(y'') \\
  &\Leftrightarrow\, A_i(y')+\sum_{j\in\{1,\dots,n\}\setminus\{i\}}A_j(y_j)
     \leq_S A_i(y'')+\sum_{j\in\{1,\dots,n\}\setminus\{i\}}A_j(y_j) \\
  &\Leftrightarrow\, 
  \Omega_1(y_1,\dots,y_{i-1},y',y_{i+1},\dots,y_n)\leq_S \Omega_1(y_1,\dots,y_{i-1},y'',y_{i+1},\dots,y_n).
}
Finally we show that $\ri(X)\subseteq \intr_{\omega_1}(X)$. Let $p\in\ri(X)$ be fixed and
define the two-variable operation $\omega^*:G^2\to G$ as $\omega^*(x,y):=a_1(x)+y-a_1(y)$. Then the identity
\Eq{*}{
  \omega^*(x,y):=a_1(x)+(a_2+\cdots+a_n)(y)=\omega_1(x,y,\dots,y) \qquad(x,y\in G)
}
and the $\omega_1$-convexity of $X$ yield that $X$ is $\omega^*$-convex. Applying \lem{ri} it follows that 
$\ri(X)\subseteq\intr_{\omega^*}(X)$, and hence $p\in\intr_{\omega^*}(X)$. By definition, this means that 
$\ext_{\omega^*}(\{p\})=X$. Now let $E\subseteq X$ be an $\omega_1$-extreme set containing $\{p\}$. We are going to 
verify that $E$ is also $\omega^*$-extreme. Indeed, if $(x,y)\in(\omega^*)^{-1}(E)$, then $\omega^*(x,y)\in E$, 
which is equivalent to $\omega_1(x,y,\dots,y)\in E$. This inclusion, by the $\omega_1$-extremality of $E$, shows that 
$(x,y,\dots,y)\in E^n$. Therefore, $(x,y)\in E^2$, which finally proves $(\omega^*)^{-1}(E)\subseteq 
E^2$, i.e., the $\omega^*$-extremality of $E$. On the other hand, we have that $\ext_{\omega^*}(\{p\})=X$, therefore 
$E=X$. Consequently, $\ext_{\omega_1}(\{p\})=X$ and thus we get that $p\in\intr_{\omega_1}(X)$.

Now we are in the position to apply the \cor{MT2}, that is all the conditions of this result are satisfied. 
Therefore, if $f:X\to Y$ is a solution of the functional inequality \eq{MT2c}, then it also fulfills the convexity 
type inequality 
\Eq{*}{
 f\big(\omega_1(x_1,\dots,x_n)\big)
    \leq_S \Omega_1\big(f(x_1),\dots,f(x_n)\big), \qquad (x_1,\dots,x_n\in X).
}
By the conclusion of \cor{MT2}, then there exists a function $g:X\to Y$ such that $g\leq_Sf$, $g(p)=f(p)$, and 
\Eq{*}{
 g\big(\omega_1(x_1,\dots,x_n)\big)
    = \Omega_1\big(g(x_1),\dots,g(x_n)\big), \qquad (x_1,\dots,x_n\in X).
}
The latter functional equation being equivalent to \eq{MT2a}, the proof of the \thm{MT2} is completed.
\end{proof}

Now, we apply the the above theorem to the proof of a support theorem for so-called delta $(s,t)$-convex maps. This theorem 
was proved in \cite{Olb15d} by the first author.  The concept of delta $(s,t)$-convex maps generalizes the concept of 
delta-convex maps which was introduced by L.\ Vesel\'y and L.\ Zaj\'{\i}\v{c}ek \cite{VesZaj89} in the following manner:
Given to real normed spaces $X$, $Y$ and a nonempty open and convex subset $D\subseteq X$, a map $F:D\to Y$ 
is said to be a \textit{delta-convex} if there exists a continuous and convex functional $f:D\to \R$ such that 
$f+y^{\star}\circ F$ is continuous and convex for any member $y^{\star}$ of the dual space of $Y$ with 
$\|y^{\star}\|=1$. If this is the case, then we say that $F$ is a \textit{delta-convex mapping with a control function 
$f$}.

It turns out that a continuous map $F:D\to Y$ is a delta-convex controlled by a continuous function $f:D\to \R$ if and 
only if the functional inequality 
\Eq{*}{
\Big\|\frac{F(x)+F(y)}{2}-F\Big(\frac{x+y}{2} \Big)\Big\|\leq \frac{f(x)+f(y)}{2}-f\Big(\frac{x+y}{2}\Big),
}
is satisfied for all $x, y\in D$.
The above functional inequality may obviously be investigated without any regularity assumptions upon $F$ and $f$ which 
additionally considerably enlarges the class of solutions. Note that, delta-convex mappings have nice properties (see 
\cite{VesZaj89}) and this notion seems to be the most natural generalization of functions which are representable as a 
difference of two convex functions. In \cite{Olb15d}, the first author generalized the concept of delta-convexity in the 
following manner: Given two numbers $s, t\in (0,1)$, a convex subset $D$ of a vector space $X$ and a Banach 
space $Y$ we say that a map $F:D\to Y$ is \emph{delta $(s,t)$-convex with a control function $f:D\to \R$}, if the 
inequality
\Eq{*}{
 \|tF(x)+(1-t)F(y)-&F(sx+(1-s)y)\|\\
 &\leq tf(x)+(1-t)f(y)-f(sx+(1-s)y),
}
holds for all $x, y\in D$.

Observe that, by defining the map $\bar{F}:D\to Y\times \R$ via the formula 
\Eq{bar}{
  \bar{F}(x):=(F(x),f(x)),\qquad (x\in D),
} 
we can rewrite the above inequality in the form
\Eq{*}{
  \bar{F}(sx+(1-s)y)\leq_{\K_1} t\bar{F}(x)+(1-t)\bar{F}(y), \qquad (x, y\in D),
}
where $\K_1:=\{(x,t)\in Y\times \R \mid \|x\|\leq t \}$ is the Lorenz cone. 

In order to formulate the main result from \cite{Olb15d}, let us recall that a map $A:D\to Y$ is said to be 
\emph{$(s,t)$-affine} if it satisfies the following functional equation 
\Eq{*}{
  A(sx+(1-s)y)=tA(x)+(1-t)A(y),\qquad (x, y\in D).
}

\Thm{st}{Let $D$ be a convex and algebraically open subset of a vector space $X$, let $Y$ be a Banach space and let 
$F:D\to Y$ be a delta $(s,t)$-convex map with a control function $f:D \to \R$. Then, for any point $y\in D$, 
there exist $(s,t)$-affine maps $A_{y}:D\to Y$ and $a_{y}:D\to\R$ such that $A_{y}(y)=F(y)$, $a_{y}(y)=f(y)$, and
\Eq{*}{
\|F(x)-A_{y}(x)\|\leq f(x)-a_{y}(x),\qquad (x\in D).
}}

\begin{proof} 
Put $\bar{Y}:=Y\times \R$ and define the map $\bar{F}:D\to\bar{Y}$ by \eq{bar}. Consider the vector ordering  
generated by the Lorenz cone $\K_1,$ which is closed, convex and sharp and consider two families of additive 
maps $a_1,a_2:X\to X$ and $A_1,A_2:\bar{Y}\to\bar{Y}$ defined by the formulas 
\Eq{*}{
 a_1(x)&:=sx,&\quad a_2(x)&:=(1-s)x,&&\qquad (x\in X); \\
 A_1(\bar{y})&:=t\bar{y},&\quad A_2(\bar{y})&:=(1-t)\bar{y},&&\qquad (\bar{y}\in \bar{Y}).
}
It is easy to see that these additive maps are commuting, moreover, 
\Eq{*}{
  a_1(x)+a_2(x)&=sx+(1-s)x=x=\id_X(x),&&\qquad (x\in X);\\ 
  A_1(\bar{y})+A_2(\bar{y})&=t\bar{y}+(1-t)\bar{y}=\bar{y}=\id_{\bar{Y}}(\bar{y}),&&\qquad (\bar{y}\in \bar{Y}).
}
Obviously, $a_1(D)+a_2(D)=sD+(1-s)D\subseteq D$ by the convexity of $D$. The operations $A_1$ 
and $A_2$ are also bijective with conditions $A_i(\K_1)=\K_1$ for $i\in \{1,2\}$. Finally, it remains to apply the 
\thm{MT2} to the inequality 
\Eq{*}{
\bar{F}(a_1(x)+a_2(y))\leq_{\K_1} A_1(\bar{F}(x))+A_2(\bar{F}(y)), \qquad (x, y\in D).
}
\end{proof}

\section*{References}


\end{document}